\documentclass[oneside]{amsart}
\usepackage[T1]{fontenc}
\usepackage[latin1]{inputenc}
\usepackage{graphicx}
\usepackage{amssymb}

\makeatletter
 \theoremstyle{plain}
\newtheorem{thm}{Theorem}[section]
  \theoremstyle{definition}
  \newtheorem{defn}[thm]{Definition}
  \theoremstyle{plain}
  \newtheorem{lem}[thm]{Lemma}
  \theoremstyle{plain}
  \newtheorem{cor}[thm]{Corollary}
  \theoremstyle{plain}
  \newtheorem{prop}[thm]{Proposition}

\newcommand{\lk}{\operatorname{lk}}

\makeatother
\begin{document}

\title{Shelling the Coset Poset}

\author{Russ Woodroofe}

\email{rsw9@cornell.edu}

\begin{abstract}
It is shown that the coset lattice of a finite group has shellable
order complex if and only if the group is complemented. Furthermore,
the coset lattice is shown to have a Cohen-Macaulay order complex
in exactly the same conditions. The group theoretical tools used are
relatively elementary, and avoid the classification of finite simple
groups and of minimal finite simple groups.
\end{abstract}
\maketitle

\section{Introduction}

We start by recalling the definition of a shelling. All posets, lattices,
simplicial complexes, and groups in this paper are finite.

\begin{defn}
If $\Delta$ is a simplicial complex, then a \emph{shelling} of $\Delta$
is an ordering $F_{1},\dots,F_{n}$ of the facets (maximal faces)
of $\Delta$ such that $F_{k}\cap\bigl(\cup_{i=1}^{k-1}F_{i}\bigr)$
is a nonempty union of facets of $F_{k}$ for $k\geq2$. If $\Delta$
has a shelling, we say it is \emph{shellable}. 
\end{defn}
We will use this definition in the context of a poset $P$ by recalling
the \emph{order complex} $|P|$ to be the simplicial complex with
vertex set $P$ and faces chains in $P$. We say that a poset $P$
is shellable if $|P|$ is. Recall also that $P$ is \emph{graded},
and $\vert P\vert$ is \emph{pure}, if all maximal chains in $P$
have the same length.

\smallskip{}
The idea of a shelling (and the property of shellability) were first
formally introduced by Bruggesser and Mani in \cite{Bruggesser/Mani:1971},
though similar ideas had been assumed implicitly since the beginning
of the 20th century. See Chapter 8 of \cite{Ziegler:1995} for a development
of some of the history and basic results on shellability. Since its
introduction, it has been studied extensively by combinatorialists.
Particularly, in the 1980's and 90's Björner and Wachs wrote several
papers \cite{Bjorner:1980,Bjorner/Wachs:1983,Bjorner/Wachs:1996,Bjorner/Wachs:1997}
developing the theory of shellability for posets. Especially important
to those interested in group theory are \cite{Bjorner/Wachs:1996}
and \cite{Bjorner/Wachs:1997}, as they extend the older definition
of shellability (which only applied to graded posets) to apply to
any poset. This extension makes Theorem \ref{thm:Shareshian} much
more interesting!

We henceforth assume that a reader has seen the basic definitions
and results of, say, \cite{Bjorner/Wachs:1996}, although we try to
state clearly what we are using.

Recall that the subgroup lattice (denoted $L(G)$) is the lattice
of all subgroups of a group $G$. Shellings of subgroup lattices have
been of interest for quite some time. In fact, one of the main results
of Björner's first paper on shelling posets \cite{Bjorner:1980} was
to show that supersolvable groups have shellable subgroup lattices.
(Recall a supersolvable group is a group having chief series with
every factor of prime order.) As mentioned before, at that time, shellability
was a property that applied only to graded posets. Under this definition,
Björner had the shellable subgroup lattices completely characterized,
if we recall the following theorem of Iwasawa:

\begin{thm}
\label{thm:Iwasawa}\emph{(Iwasawa \cite{Iwasawa:1941})} Let $G$
be a finite group. Then $L(G)$ is graded if and only if $G$ is supersolvable.
\end{thm}
Of course, when Björner and Wachs updated the definition of a shelling
to allow non-graded posets in \cite{Bjorner/Wachs:1996,Bjorner/Wachs:1997},
shellable subgroup lattices were no longer characterized. This gap
was soon filled by Shareshian:

\begin{thm}
\emph{\label{thm:Shareshian}(Shareshian \cite{Shareshian:2001})}
Let $G$ be a finite group. Then the subgroup lattice $L(G)$ is shellable
if and only if $G$ is solvable.
\end{thm}
A nice summary article on shellability and group theory was written
by Welker in \cite{Welker:1994}. This article is now somewhat out
of date, and it has some errors, but it is very useful as an overview
of the topic. The reader should be warned, however, that at the time
it was written shellability was still considered to apply only to
graded posets.

\bigskip{}
Shareshian's result is surprising and pretty, and it would be nice
to find something similar for other lattices on groups. In this paper,
we consider cosets. The \emph{coset poset} $\mathfrak{C}(G)$ (poetically
named by K.~Brown in \cite{Brown:2000}) is the poset of all cosets
of proper subgroups of $G$, ordered by inclusion. The \emph{coset
lattice} $\hat{\mathfrak{C}}(G)$ is $\mathfrak{C}(G)\cup\{ G,\emptyset\}$,
likewise ordered by inclusion. The meet operation is intersection,
while $H_{1}x_{1}\vee H_{2}x_{2}=\langle H_{1},H_{2},x_{1}x_{2}^{-1}\rangle x_{2}$.
Clearly, $\mathfrak{C}(G)$ is shellable if and only if $\hat{\mathfrak{C}}(G)$
is, so we study the two interchangeably. If $\mathfrak{C}(G)$ is
shellable, we will call $G$ \emph{coset-shellable}.

The history of the coset poset is discussed in the last chapter of
\cite{Schmidt:1994}. Most results proved have been either negative
results, or else so similar to the situation in the subgroup lattice
as to be uninteresting. More recently, K.~Brown rediscovered the
coset poset, and studied its homotopy type while proving some divisibility
results on the so-called probabilistic zeta function \cite{Brown:2000}.
In particular, he showed that if $G$ is a solvable group, then $|\mathfrak{C}(G)|$
has the homotopy type of a bouquet of spheres, all of the same dimension. 

In Section \ref{sec:NonshellableCPs} we show that there are finite
groups $G$ which have a shellable subgroup lattice, but a non-shellable
coset lattice. In particular, we show that for $\mathfrak{C}(G)$
to be shellable, $G$ must be supersolvable, and every Sylow subgroup
of $G$ must be elementary abelian. Our main tool is the above-mentioned
result of Brown, together with the fact that a pure shellable complex
is homotopic to a bouquet of top-dimensional spheres. In Section \ref{sec:LinearAlgebra}
we use linear algebra to construct an invariant on subgroups of such
groups. Finally, in Section \ref{sec:ShellingCPs} we use this invariant
to construct a so-called EL-shelling, and to finish the proof of our
main theorem:

\begin{thm}
\emph{(Main Theorem)} If $G$ is a finite group, then $\mathfrak{C}(G)$
is shellable if and only if $G$ is supersolvable with all Sylow subgroups
elementary abelian.
\end{thm}
Our theorem is even more interesting when we connect it with a paper
of P.~Hall \cite{PHall:1937}. We recall that a group $G$ is \emph{complemented}
if for every subgroup $H\subseteq G$, there is a \emph{complement}
$K$ which satisfies i) $K\cap H=1$ and ii) $HK=KH=G$. Hall proved
the equivalence of the first three properties in the following restatement
of our theorem:

\begin{thm}
\emph{\label{thm:MainThmRestated}(Restatement of Main Theorem)} If
$G$ is a finite group, then the following are equivalent:
\begin{enumerate}
\item $G$ is supersolvable with all Sylow subgroups elementary abelian,
\item $G$ is complemented,
\item $G$ is a subgroup of the direct product of groups of square free
order,
\item $G$ is coset-shellable,
\item $\mathfrak{C}(G)$ is homotopy Cohen-Macaulay,
\item $\mathfrak{C}(G)$ is sequentially Cohen-Macaulay over some field,
and
\item $\mathfrak{C}(G)$ is Cohen-Macaulay over some field.
\end{enumerate}
\end{thm}
Parts (5), (6) and (7) are discussed in Section \ref{sub:Cohen-Macaulay-coset},
where we define the three used versions of the Cohen-Macaulay property,
and give further references.

Notice that, contrary to the situation of Shareshian's Theorem, $\mathfrak{C}(G)$
is shellable if and only if it is pure and shellable. Thus, non-pure
shellability only comes in in the negative direction of our proof.
Of course, now that it has been defined, one can't ignore it!

Complemented groups have also been called \emph{completely factorizable}
groups, and have been studied by other people, see for example \cite{Ballester-Bolinches/Xiuyun:1999},
or \cite{Malanina/Hlebutina/Sevcov:1972}. Ramras has further examined
the homotopy type of the coset poset in \cite{Ramras:2005}.

I would like to thank my advisor, Ken Brown, for many discussions
and ideas; and Keith Dennis for helping with some of my group theory
questions. Thanks also to Yoav Segev, Anders Björner, and the anonymous
referees for their helpful comments. Though no computations appear
in this paper, GAP \cite{GAP4} and the XGAP package often helped
in searching for examples that furthered understanding.

\section{Coset posets that are not shellable\label{sec:NonshellableCPs}}

\subsection{$p$-groups}

It is often easier to show that something is not shellable, than to
show that it is. So we start our search for shellings of the coset
lattice by finding groups for which $\mathfrak{C}(G)$ is certainly
\emph{not} shellable. The following lemma will be very useful in this
endeavor.

\begin{lem}
\label{lem:ShellableIntervals} If $P$ is a shellable poset, then
every interval in $P$ is also shellable. (Thus, if $G$ is coset-shellable,
then so is every subgroup $H\subseteq G$.)
\end{lem}
\begin{proof}
Since every interval in a poset $P$ is a {}``link'' (for more information,
see the beginning of Section \ref{sub:Cohen-Macaulay-coset}), the
first part follows immediately from Proposition 10.14 in \cite{Bjorner/Wachs:1997}.

For the second part, we note that the interval $[\emptyset,H]$ in
$\hat{\mathfrak{C}}(G)$ is isomorphic to $\hat{\mathfrak{C}}(H)$.
\end{proof}
\begin{cor}
\label{cor:CShellableIsSolvable}If $G$ is a finite coset-shellable
group, then $G$ is solvable.
\end{cor}
\begin{proof}
Note that the interval $[1,G]$ in $\hat{\mathfrak{C}}(G)$ is isomorphic
to the subgroup lattice of $G$. Apply Lemma \ref{lem:ShellableIntervals}
and Shareshian's Theorem (Theorem \ref{thm:Shareshian}).
\end{proof}
A proof of Corollary \ref{cor:CShellableIsSolvable} that does not
rely on Shareshian's Theorem will also be given, in Section \ref{sub:NonSupersolvableGroups}.

\begin{figure}[htbp]
\begin{centering}\includegraphics{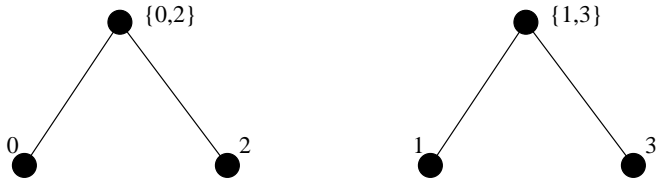}\par\end{centering}

\caption{The coset poset of $\mathbb{Z}_{4}$\label{fig:z4} is not connected,
so not shellable.}
\end{figure}

At first glance, one might hope that perhaps all solvable groups have
a shellable coset poset. Soon enough, however, one considers the coset-poset
of $\mathbb{Z}_{4}$, pictured in figure \ref{fig:z4}. We see that
$\mathfrak{C}(\mathbb{Z}_{4})$ is not even connected, and connectivity
is an easy consequence of the definition of shellability as long as
all facets have dimension at least 1. 

A similar situation holds for an arbitrary prime $p$: $\mathbb{Z}_{p^{2}}$
has only one nontrivial proper subgroup, so $\mathfrak{C}(\mathbb{Z}_{p^{2}})$
falls into $p$ connected components, and in particular is not shellable.
Hence, no group $G$ with a subgroup isomorphic with $\mathbb{Z}_{p^{2}}$
can be coset-shellable. Can we eliminate any other $p$-groups from
the possibility of coset-shellability? In fact we can. We recall the
following theorem of K.~Brown:

\begin{thm}
\emph{\label{thm:HomotopyTypeForSolvable}(Proposition 11 from \cite{Brown:2000})}
Let $G$ be a finite solvable group with a chief series $1=N_{0}\subseteq N_{1}\subseteq\dots\subseteq N_{k}=G$.
Then $\mathfrak{C}(G)$ has the homotopy type of a bouquet of $(d-1)$-spheres,
where $d$ is the number of indices $i=1,\dots,k$ such that $N_{i}/N_{i-1}$
has a complement in $G/N_{i-1}$.
\end{thm}
It follows from the proof in \emph{\cite{Brown:2000}} that for $G\neq1$
the number of spheres is, in fact, nonzero. (It is \[
\prod_{i=1}^{k}(1-c_{i}\vert N_{i}/N_{i-1}\vert)\]
where $c_{i}$ is the number of complements $N_{i}/N_{i-1}$ has in
$G/N_{i-1}$; also related is \cite[Corollary 3]{Brown:2000}.)

\begin{prop}
\label{pro:CShellableSubgroupsAreEAbel}~\\
1) If $H$ is a finite $p$-group which is coset-shellable, then $H$
is elementary abelian.\\
2) If $G$ is a finite group which is coset-shellable, then all Sylow
subgroups of $G$ are elementary abelian.
\end{prop}
\begin{proof}
1) If $H$ is a finite $p$-group, then $L(H)$ is graded (by, for
example, Iwasawa's Theorem, Theorem \ref{thm:Iwasawa}), hence $\mathfrak{C}(H)$
is also graded. But it is well known (see for example \cite{Bjorner/Wachs:1996})
that a graded, shellable poset $P$ has homotopy type of a bouquet
of $r$-spheres, where $r$ is the length of a maximal chain in $P$.
By the above theorem, we see that every chief factor of $H$ must
be complemented, and hence $H$ has trivial Frattini subgroup $\Phi(H)$
(otherwise any minimal normal subgroup contained in $\Phi(H)$ is
an uncomplemented chief factor).

But for a finite $p$-group, $\Phi(H)=H'H^{p}$ (see for example \cite[5.3.2]{Robinson:1996}),
so $H$ is abelian of exponent $p$, that is, elementary abelian.

2) Apply Lemma \ref{lem:ShellableIntervals} to the interval $[\emptyset,H]$
in $\hat{\mathfrak{C}}(G)$, where $H$ is a Sylow subgroup of $G$. 
\end{proof}

\subsection{Non-supersolvable groups\label{sub:NonSupersolvableGroups}}

We now have that for a finite group $G$ to be coset-shellable, $G$
must be solvable with elementary abelian Sylow subgroups. A little
more holds: $G$ must in fact be supersolvable. To prove this, it
suffices by Lemma \ref{lem:ShellableIntervals} and the discussion
of the previous section to restrict ourselves to groups $G$ such
that: 

\begin{enumerate}
\item $G$ is not supersolvable,
\item All proper subgroups of $G$ are supersolvable, and
\item All Sylow subgroups of $G$ are elementary abelian.
\end{enumerate}
A closely related idea is that of minimal non-complemented groups,
which are non-complemented groups with every proper subgroup complemented.
A complete characterization of such groups is given in \cite{Malanina/Hlebutina/Sevcov:1972},
although we do not use their characterization.

In light of Shareshian's Theorem and Corollary \ref{cor:CShellableIsSolvable},
it might seem at first glance that a stronger version of Condition
1 would be to require $G$ to be solvable but not supersolvable. The
following result of Doerk, however, shows that this would be redundant.

\begin{lem}
\label{lem:MinimalNonsupersolvableAreSolvable}\emph{(Doerk, \cite[Hilfssatz C]{Doerk:1966})}
If every maximal subgroup of $G$ is supersolvable, then $G$ is solvable.
\end{lem}
We notice also that this frees our characterization of groups that
are not coset-shellable from Shareshian's Theorem, which relies on
Thompson's classification of minimal finite simple groups. We will
also need this for Section \ref{sub:Cohen-Macaulay-coset}.

For any normal subgroup $N$ in $G$, let $q:G\rightarrow G/N$ be
the quotient map. Then we take $\mathfrak{C}_{0}(G)$ to be the subposet
of $\mathfrak{C}(G)$ of all $Hx$ such that $q(Hx)\neq G/N$. Thus,
$\mathfrak{C}_{0}(G)$ is obtained from $\mathfrak{C}(G)$ by removing
cosets $Kx$ when $KN=G$. We will use the following proposition to
show that, for $G$ satisfying Conditions 1--3, $\mathfrak{C}(G)$
has the wrong homotopy type to be shellable.

\begin{thm}
\emph{\label{thm:BrownC_0Homotopy}(K.~Brown \cite{Brown:2000},
Proposition 8 and following discussion)} The quotient map $q:G\rightarrow G/N$
induces a homotopy equivalence $\mathfrak{C}_{0}(G)\rightarrow\mathfrak{C}(G/N)$.
\end{thm}
The following lemma from group theory will be useful.

\begin{lem}
\label{lem:ComplementsOfChiefFactors}Let $G$ be a solvable group,
with $H$ a proper subgroup. Then
\begin{enumerate}
\item If $N$ is an abelian normal subgroup of $G$ with $NH=G$, then H
is maximal in $G$ if and only if $N/N\cap H$ is a chief factor for
$G$.
\item $H$ is maximal if and only if it is a complement to a chief factor
$N_{i+1}/N_{i}$, i.e., if and only if $HN_{i+1}=G$ and $H\cap N_{i+1}=N_{i}$.
\end{enumerate}
\end{lem}
Part 1 may be found in \cite[Theorem 5.4.2]{Robinson:1996}. Part
2 follows from part 1 by taking $N_{i+1}$ to be minimal such that
$HN_{i+1}=G$.

We use Lemma \ref{lem:ComplementsOfChiefFactors} in proving the following:

\begin{lem}
\label{lem:C_0IsPureSkeleton}Let $G$ be a group satisfying Conditions
1--3 above. Let $n$ be the length of a longest chain in $\mathfrak{C}(G)$.
Then $G$ has a minimal normal subgroup $N$, of non-prime order,
over which $\vert\mathfrak{C}_{0}(G)\vert$ is the subcomplex of $\vert\mathfrak{C}(G)\vert$
generated by all chains of length $n$.
\end{lem}
\begin{proof}
Our proof goes in five steps:

1) Every chief factor $N_{i+1}/N_{i}$ of $G$ is complemented in
$G/N_{i}$.

We apply a theorem of Gaschütz, proved in \cite[Theorem 3.3.2]{Kurzweil/Stellmacher:2004},
which says that a normal abelian $p$-subgroup $N$ has a complement
in $G$ if and only if $N$ has a complement in a Sylow $p$-subgroup
$P$ containing $N$. Let $N_{i+1}/N_{i}$ be a chief factor of $G$.
Then since $G$ has all Sylow subgroups elementary abelian, $G/N_{i}$
has all Sylow subgroups elementary abelian. But an elementary abelian
group is a complemented group (see Theorem \ref{thm:MainThmRestated}),
hence $N_{i+1}/N_{i}$ has a complement in any Sylow subgroup $P/N_{i}$
containing it, and so by Gaschütz we get that $N_{i+1}/N_{i}$ has
a complement in $G/N_{i}$.

2) A chief factor $N_{i+1}/N_{i}$ is of non-prime order only if $N_{i}=1$.

Suppose otherwise, that $N_{i+1}/N_{i}$ is of non-prime order with
$N_{i}\neq1$, so that $G/N_{i}$ is solvable but not supersolvable.
Then $N_{i}/N_{i-1}$ has a complement by part (1), so there is a
group $K$ with $G/N_{i}\cong K/N_{i-1}$. Since all subgroups of
$G$ are supersolvable, we see that $G/N_{i}$ is supersolvable, a
contradiction.

3) There exists a minimal normal subgroup $N\subseteq G$ of non-prime
order, and $N$ is a complement to any maximal subgroup $K$ of non-prime
index.

Since $G$ is not supersolvable, there is some factor of non-prime
order in any chief series of $G$, and by part (2) it must be of the
form $N_{1}/1$. We notice that in the situation of the second part
of Lemma \ref{lem:ComplementsOfChiefFactors}, we have $[G:H]=[N_{i+1}:N_{i}]$.
Since $N_{1}/1$ is the only factor of non-prime order, we get the
desired result.

4) A maximal chain $C$ has length less than $n$ if and only if the
top element of $C$ is a coset $Kx$ of some complement $K$ of $N$. 

Suppose $Hx$ is the top element of a chain $C$ in $\mathfrak{C}(G)$.
By Iwasawa's Theorem and the supersolvability of $H$, we have that
every maximal chain in the interval $(\emptyset,Hx]$ has length $n+1-a$,
where $a$ is the number of primes with multiplicity dividing $[G:H]$.
We see that $C$ has length $n$ if and only if $a=1$, so $C$ is
of length less than $n$ if and only if $H$ is of non-prime index
in $G$ if and only if $H$ is a complement of $N$.

5) In the situation of part (4), $C\setminus\{ Kx\}$ can be extended
to a chain of length $n$.

Let $K_{1}x$ be the coset immediately under $Kx$ in $C$. Then since
$K$ is supersolvable, $[K:K_{1}]$ is a prime. Then $[G:NK_{1}]=[K:K_{1}]$
is also a prime. Moreover, $NK_{1}$ is supersolvable, so if $\vert N\vert=p^{a}$,
then there is a chain $K_{1}=H_{0}<H_{1}<\dots<H_{a}=NK_{1}$ between
$K_{1}$ and $NK_{1}$. The desired chain then follows $C$ up to
$K_{1}x$, and ends at the top with $K_{1}x<H_{1}x<\dots<H_{a}x=NK_{1}x$.

\smallskip{}

We have shown that $\vert\mathfrak{C}_{0}(G)\vert$ is obtained from
$\vert\mathfrak{C}(G)\vert$ by removing the facets of dimension less
than $n$, thus that $\vert\mathfrak{C}_{0}(G)\vert$ is the subcomplex
of $\vert\mathfrak{C}(G)\vert$ generated by all $n$-faces.
\end{proof}
We notice in passing that the argument in part (2) actually shows
that a complement $K$ of $N$ has $\operatorname{Core}_{G}K=1$,
thus that a group $G$ satisfying Conditions 1--3 is \emph{primitive}.
Such groups have highly restricted structure, see Chapter A.15 of
\cite{Doerk/Hawkes:1992} for an overview.

We relate the preceding lemma to the following result from Björner
and Wachs:

\begin{lem}
\label{lem:SkeletonsShellable}\emph{(Björner/Wachs \cite[Theorem 2.9]{Bjorner/Wachs:1996})}
If $\Delta$ is shellable, then the subcomplex generated by all faces
of dimensions between $r$ and $s$ is also shellable, for all $r\leq s$.
\end{lem}
We are now ready to prove our goal for this section.

\begin{thm}
\label{thm:NotSupersolvableNotShellable}If $G$ is not supersolvable,
then $G$ is not coset-shellable.
\end{thm}
\begin{proof}
By the preceding discussion, it suffices to consider $G$ solvable
with every subgroup a complemented group. Let $N$ be the minimal
normal subgroup of order $p^{a}$ constructed in Lemma \ref{lem:C_0IsPureSkeleton}.
Then the resulting $\vert\mathfrak{C}_{0}(G)\vert$ is the subcomplex
of $\vert\mathfrak{C}(G)\vert$ generated by the faces of dimension
$n$. Theorem \ref{thm:BrownC_0Homotopy} gives us that $\vert\mathfrak{C}_{0}(G)\vert\simeq\vert\mathfrak{C}(G/N)\vert$. 

If $\mathfrak{C}_{0}(G)$ were shellable, then $\vert\mathfrak{C}_{0}(G)\vert$
would have the homotopy type of a bouquet of $n$-spheres, as discussed
in the proof of Proposition \ref{pro:CShellableSubgroupsAreEAbel}.
But Theorem \ref{thm:HomotopyTypeForSolvable} and the comment following
give that $\vert\mathfrak{C}_{0}(G)\vert\simeq\vert\mathfrak{C}(G/N)\vert$
is homotopic to a non-empty bouquet of $(n-a)$-spheres. Thus, $\mathfrak{C}_{0}(G)$
is not shellable, and by Lemma \ref{lem:SkeletonsShellable} we see
that $\mathfrak{C}(G)$ is not shellable.
\end{proof}
We have now proved that (4) $\implies$ (1--3) in our Restatement
of the Main Theorem, Theorem \ref{thm:MainThmRestated}. The following
subsection, which the rest of the paper does not depend on, deals
with (5--7). A reader who is unfamiliar with the Cohen-Macaulay property
for simplicial complexes may, if desired, skip directly to Section
\ref{sec:LinearAlgebra}.

\subsection{Cohen-Macaulay coset lattices\label{sub:Cohen-Macaulay-coset}}

In fact, we have proven slightly more in Section \ref{sub:NonSupersolvableGroups}.
Two properties that are closely related to shellability are that of
being Cohen-Macaulay and (generalized to non-pure complexes) that
of being sequentially Cohen-Macaulay. Recall that the \emph{link}
of a face $F_{0}$ in a simplicial complex $\Delta$ is $\lk_{\Delta}F_{0}=\{ F\in\Delta:F\cup F_{0}\in\Delta,\, F\cap F_{0}=\emptyset\}$.
Links in the order complexes of posets are closely related to intervals.
More specifically, if $C$ is a maximal chain containing $x$ and
$y$, and $C'$ is $C$ with all $z$ such that $x<z<y$ removed,
then it is easy to see that $\lk_{P}C'$ is the order complex of the
interval $(x,y)$. In general, the link of a chain in a bounded poset
is the so-called {}``join'' of intervals.

Let $k$ be a field. A simplicial complex $\Delta$ is \emph{Cohen-Macaulay
over $k$} if for every face $F\in\Delta$, $\tilde{H}_{i}(\lk_{\Delta}F,k)=0$
for $i<\dim\lk_{\Delta}F$, i.e., if every link has the homology of
a wedge of top dimensional spheres. It will come as little surprise
after the preceding discussion of links in posets that one can prove
the following fact: a poset $P$ is Cohen-Macaulay if and only if
every interval $(x,y)$ in $P$ has the homological wedge of spheres
property (see \cite{Bjorner:1995} for a proof of this and further
discussion of links and joins). The complex $\Delta$ is \emph{homotopy
Cohen-Macaulay} if every such link is homotopic to (rather than merely
having the homology of) a wedge of top dimensional spheres. Since
a graded shellable poset has the homotopy type of a wedge of top dimensional
spheres, and since every interval in a shellable poset is shellable,
we see that (the order complex of) a graded shellable poset is homotopy
Cohen-Macaulay. A homotopy Cohen-Macaulay complex is Cohen-Macaulay
over any field.

There is an extension of the Cohen-Macaulay property to non-pure complexes.
The \emph{pure $i$-skeleton} of a simplicial complex $\Delta$ is
the subcomplex generated by all faces of dimension $i$. We say that
$\Delta$ is \emph{sequentially Cohen-Macaulay} if its pure $i$-skeleton
is Cohen-Macaulay for all $i$. A pure, sequentially Cohen-Macaulay
complex is obviously Cohen-Macaulay. 

A reference for background on Cohen-Macaulay complexes is \cite{Stanley:1996}.
Useful properties of sequentially Cohen-Macaulay complexes are given
in \cite{Wachs:1999}. We recall some facts presented in the latter.

\begin{lem}
\label{lem:CMproperties}Let $\Delta$ be a simplicial complex, $P$
a poset:
\begin{enumerate}
\item If $\Delta$ is shellable, then $\Delta$ is sequentially Cohen-Macaulay.
\cite[Corollary 1.6]{Wachs:1999}
\item If $P$ is sequentially Cohen-Macaulay, then all intervals in $P$
are also sequentially Cohen-Macaulay. \cite[Theorem 1.5]{Wachs:1999}
\end{enumerate}
\end{lem}
Then in the previous two sections we have actually shown

\begin{prop}
\label{pro:CMareComplemented}If $\mathfrak{C}(G)$ is sequentially
Cohen-Macaulay, then $G$ is a complemented group.
\end{prop}
\begin{proof}
The proof of Proposition \ref{pro:CShellableSubgroupsAreEAbel} shows
that if $P$ is a $p$-group, but not elementary abelian, then $\mathfrak{C}(P)$
has the homotopy type of a wedge of spheres of the wrong dimension.
Hence the homology does not vanish below the top dimension, and $\mathfrak{C}(P)$
is not (sequentially) Cohen-Macaulay. Lemma \ref{lem:CMproperties}
part 2 then gives that all Sylow subgroups of a group $G$ with $\mathfrak{C}(G)$
sequentially Cohen-Macaulay must be elementary abelian.

Similarly, in the proof of Theorem \ref{thm:NotSupersolvableNotShellable}
we show that $\mathfrak{C}_{0}(G)$ is not Cohen-Macaulay. By Lemma
\ref{lem:C_0IsPureSkeleton} we have that $\mathfrak{C}_{0}(G)$ is
the pure $n$-skeleton of $\mathfrak{C}(G)$, and then the definition
gives that $\mathfrak{C}(G)$ is not sequentially Cohen-Macaulay unless
$G$ is supersolvable.
\end{proof}
Proposition \ref{pro:CMareComplemented} and the fact that complemented
groups are supersolvable then give that (6) and (7) are equivalent,
and that both imply (1--3) in our Restatement of the Main Theorem.
Then (4) $\implies$(5) $\implies$(6) is clear from the definition
of homotopy Cohen-Macaulay, and it remains only to prove (1--3) $\implies$(4).
This will be the subject of Section \ref{sec:ShellingCPs}.

\section{Some linear algebra\label{sec:LinearAlgebra}}

We now take a brief break from shellings and homotopy type to do some
linear algebra. First, we introduce some notation. Fix a vector space
$V$ with (ordered) basis $\mathfrak{B}=\{ e_{1},\dots,e_{n}\}$,
and consider a subspace $U\subseteq V$. Let $\{ g_{1},\dots,g_{k}\}$
be a set of generators for $U$. Then we can write the coordinates
of the $g_{i}$'s as row vectors $[g_{i}]_{\mathfrak{B}}$, put these
in a matrix $\left[\begin{array}{c}
g_{1}\\
\vdots\\
g_{k}\end{array}\right]_{\mathfrak{B}}$, and reduce to reduced row echelon form $M$. Denote the set of pivot
columns for $M$ (i.e., the columns with a leading 1 in some row of
$M$) as $I_{V,\mathfrak{B}}(U)$, or just $I(U)$ if the choice of
$V$ and $\mathfrak{B}$ is clear.

\begin{lem}
$I(U)$ is an invariant for the subspace $U$ of $V$ with respect
to \textbf{$\mathfrak{B}$}.
\end{lem}
\begin{proof}
We need only show that $I(U)$ does not depend on the choice of generators
for $U$. Suppose generators $\{ h_{i}\}$ give row reduced matrix
$M_{h}$ and generators $\{ g_{i}\}$ give row reduced matrix $M_{g}$.
But then the row reduced matrix of $\{ h_{i}\}\cup\{ g_{i}\}$ must
be both $M_{h}$ and $M_{g}$ by uniqueness of reduced row echelon
form, hence $M_{g}=M_{h}$. In particular, the pivot columns are the
same.
\end{proof}
We mention some elementary properties of our invariant.

\begin{prop}
\label{prop:PropertiesOfI}Fix $V$ and $\mathfrak{B}$ as above,
and let $U_{1},U_{2}$ be subspaces of $V$. Then
\begin{enumerate}
\item $|I(U_{1})|=\dim U_{1}$.
\item If $U_{1}\subseteq U_{2}$, then $I(U_{1})\subseteq I(U_{2})$.
\end{enumerate}
\end{prop}
\begin{proof}
From a first course in linear algebra, the number of pivots of a matrix
is the dimension of the row space, and adding rows to the matrix adds
pivots, but does not change the ones we had before.
\end{proof}
We will need the following lemma in our application of $I(U)$ to
the next section. Briefly, part (2) will correspond with having a
unique lexicographically first path in intervals of $\mathfrak{C}(G)$.

\begin{lem}
\emph{\label{lem:TechnicalLemma}}Fix $V$ and $\mathfrak{B}$ as
above, and let $U_{1}\subseteq U_{2}$ be subspaces of $V$. Then
\begin{enumerate}
\item If $k$ is the largest number in $I(U_{2})\setminus I(U_{1})$, then
there is a unique subspace $W_{\uparrow}$ such that $U_{1}\subseteq W_{\uparrow}\subseteq U_{2}$
and $I(W_{\uparrow})=I(U_{1})\cup\{ k\}$.
\item If $j$ is the smallest number in $I(U_{2})\setminus I(U_{1})$, then
there is a unique subspace $W_{\downarrow}$ such that $U_{1}\subseteq W_{\downarrow}\subseteq U_{2}$
and $I(W_{\downarrow})=I(U_{2})\setminus\{ j\}$.
\end{enumerate}
\end{lem}
\begin{proof}
1) It is immediate from the definition of $I(U_{2})$ that there is
some $g\in U_{2}$ with a 1 in the $k$th coordinate, and 0's in all
preceding coordinates when written as a vector with respect to $\mathfrak{B}$.
Suppose $g_{1}$ and $g_{2}$ both have this property. Then $g_{1}-g_{2}$
is 0 in all coordinates up to and including $k$, hence $g_{1}-g_{2}\in U_{1}$.
It follows that the desired $W_{\uparrow}=\langle U_{1},g\rangle$
is unique.

2) First, such a subspace exists. Suppose $U_{2}=\langle g_{1},\dots,g_{n}\rangle$,
where the $g_{i}$'s are row reduced as in the definition of $I(U)$.
Reorder so that $g_{1},\dots,g_{l}$ are the generators (rows) with
pivots in $I(U_{2})\setminus I(U_{1})$, ordered from least to greatest
(where $l=\vert I(U_{2})\setminus I(U_{1})\vert$). Then $U_{2}=\langle U_{1},g_{1},\dots,g_{l}\rangle$,
and $W_{\downarrow}=\langle U_{1},g_{2},\dots,g_{l}\rangle$ is a
space with the desired properties.

Suppose $W$ is another such space. Represent $W=\langle U_{1},h_{2},\dots,h_{l}\rangle$
in the same way as we did for $U_{2}$ in the preceding paragraph.
Let $W_{0}=\langle g_{2},\dots,g_{l},h_{2},\dots,h_{l}\rangle$. Then
the $h_{i}$'s and $g_{i}$'s are all zero in coordinates up to and
including $j$, so $j\notin I(W_{0})$. Also, $W_{0}\subseteq U_{2}$
so $I(W_{0})\subseteq I(U_{2})$. But the $g_{i}$'s and $h_{i}$'s
were row reduced with respect to $U_{1}$, so are zero in all pivots
of $U_{1}$, so $I(U_{1})\cap I(W_{0})=\emptyset$. To summarize,
$I(W_{0})\subseteq I(U_{2})\setminus\left(\{ j\}\cup I(U_{1})\right)$.
But since $W_{0}$ is at least $l-1$ dimensional (as the $g_{i}$'s
are linearly independent), we get that this is actually an equality.
Thus, $\langle g_{2},\dots,g_{l}\rangle=\langle h_{2},\dots,h_{l}\rangle=W_{0}$,
and we have that $W=W_{\downarrow}$.
\end{proof}
This ends our excursion into linear algebra. We are now ready to apply
the results of this section.

\section{Shelling the Coset Poset\label{sec:ShellingCPs}}

To show that the coset poset of a finite complemented group $G$ is
shellable, we actually exhibit a co$EL$-labeling. First, let us recall
the definition of an $EL$-labeling. 

A \emph{cover relation} is a pair $x\leftarrow y$ in a poset $P$
such that $x\lvertneqq y$ and such that there is no $z$ with $x\lvertneqq z\lvertneqq y$.
In this situation, we say that $y$ \emph{covers} $x$. We recall
that the usual picture one draws of a poset $P$ is the \emph{Hasse
diagram}, where we arrange vertices corresponding with the elements
of $P$ such that $x$ is below $y$ if $x<y$, and draw an edge between
$x$ and $y$ if $x\leftarrow y$. We say that a poset is \emph{bounded}
if it has a unique top and bottom, that is, unique upper and lower
bounds.

Let $\lambda$ be a labeling of the cover relations (equivalently,
of the edges of the Hasse diagram) of $P$ with elements of some poset
$L$ -- for us, $L$ will always be the integers. Then $\lambda$
is an \emph{$EL$-labeling} if for every interval $[x,y]$ of $P$
we have i) there is a unique (strictly) increasing maximal chain on
$[x,y]$, and ii) this maximal chain is first among maximal chains
on $[x,y]$ with respect to the lexicographic ordering. If $\lambda$
is an $EL$-labeling of the dual of $P$, then we say $\lambda$ is
a \emph{co$EL$-labeling}. 

Björner first introduced $EL$-labelings in \cite{Bjorner:1980},
and showed that if a bounded poset $P$ has an $EL$-labeling, then
$P$ is shellable. For this reason, posets with an $EL$-labeling
(or co$EL$-labeling) are often called \emph{$EL$-shellable} (or
\emph{co$EL$-shellable}). As we mentioned before, we will use the
invariants $I(U)$ discussed in the previous section to construct
a co$EL$-labeling of $\hat{\mathfrak{C}}(G)$.

\smallskip{}
Let $G=G_{1}\times\dots\times G_{r}$ be the direct product of square
free groups $\{ G_{i}\}$, and identify each group $G_{i}$ with its
inclusion in $G$. Fix $p$ a prime. Let $H$ be a subgroup of $G$,
with $H^{*}$ a Sylow $p$-subgroup of $H$. Let $G^{*}$ be a Sylow
$p$-subgroup of $G$, with $H^{*}$ contained in $G^{*}$. By the
normality of $G_{i}$ and an order argument, we get that $G^{*}\cap G_{i}$
is either isomorphic to $\mathbb{Z}_{p}$ or $1$ (depending on whether
$p\,|\,\left|G_{i}\right|$). Let $e_{i}$ be a generator of $G^{*}\cap G_{i}$
when this intersection is nontrivial. Let $\mathfrak{B}$ be an ordered
basis of such generators $e_{i}$, taken from each nontrivial $G^{*}\cap G_{i}$.
Think of the elementary abelian subgroup $G^{*}$ as a vector space
over $\mathbb{Z}_{p}$, and define $I^{p}(H)$ to be $I_{G^{*},\mathfrak{B}}(H^{*})$.

\begin{lem}
\label{lem:I^p(H)-iswelldefined}$I^{p}(H)$ is well-defined.
\end{lem}
\begin{proof}
We need to check that $I^{p}(H)$ is independent of the choice of
$H^{*}$, $G^{*}$, and $\mathfrak{B}$. Recall that $I_{G^{*},\mathfrak{B}}(H^{*})$,
as defined in Section \ref{sec:LinearAlgebra}, is the set of pivots
of the matrix with rows generating $H^{*}$. 

Notice that an element $g=e_{1}^{\alpha_{1}}\cdots e_{r}^{\alpha_{r}}\in G^{*}$
has $\alpha_{j}=0$ for $j<i$ if and only if $g\in G_{i}G_{i+1}\cdots G_{r}$.
Now the pivot associated with $e_{i}$ is in $I_{G^{*},\mathfrak{B}}(H^{*})$
if and only if there is an element in the matrix of row-reduced generators
for $H^{*}$ with first non-zero position $i$, that is, of the form
$h=e_{i}^{\alpha_{i}}\cdots e_{r}^{\alpha_{r}}$ with $\alpha_{i}\neq0$.
We see that this happens if and only if $h\in G_{i}\cdots G_{r}\setminus G_{i+1}\cdots G_{r}$.
We also notice that a set of generators is (weakly) row-reduced if
and only if the first non-zero positions are strictly increasing.
Thus, the set of pivots for $H^{*}$ is determined by the subgroups
of the form $G_{j}\cdots G_{r}$, and as $G_{j}\cdots G_{r}$ does
not depend on $G^{*}$ or $\mathfrak{B}$, it follows immediately
that $I^{p}(H)$ is independent of them.

It remains to check that $I^{p}(H)$ is independent of the choice
of $H^{*}$. Any alternate choice differs only by conjugation by some
element $x\in H$. But a set of generators $\{ h_{1},\dots,h_{k}\}$
for $H^{*}$ has the same matrix representation with respect to $\mathfrak{B}$
as their conjugates $\{ x^{-1}h_{1}x,\dots,x^{-1}h_{k}x\}$ has with
respect to $x^{-1}\mathfrak{B}x$. In particular, the set of pivots
is unchanged.
\end{proof}
We need a couple more lemmas.

\begin{lem}
\label{lem:MaxlSubgpsIntersect}Let $M$ be a maximal subgroup of
a supersolvable group $G$. If $G=HM$, then $Hx\cap M$ is a maximal
coset of $Hx$.
\end{lem}
\begin{proof}
Since $G=HM$, we can write $Hx=Hm$ for some $m\in M$. So $Hx\cap M=(H\cap M)m\neq\emptyset$.
Also, $|G|=|HM|=\frac{|H|\,|M|}{|H\cap M|}=[H:H\cap M]\cdot\frac{|G|}{[G:M]}$.
Since $[G:M]$ is prime, it follows that $[H:H\cap M]=[G:M]$ is also
prime, hence that $H\cap M$ is maximal in $H$.
\end{proof}
The following lemma (due to G.~Zappa) is proved, for example, in
\cite[5.4.8]{Robinson:1996}.

\begin{lem}
\label{lem:NormalPSubgpsOfSS}Let $G$ be a finite supersolvable group.
Then $G$ has a chief series $1=N_{0}\subseteq\dots\subseteq N_{k}=G$
with $[N_{1}:N_{0}]\geq[N_{2}:N_{1}]\geq\dots\geq[N_{k}:N_{k-1}]$.
\end{lem}
In particular, if $p$ is the largest prime dividing $\vert G\vert$
and $q$ is the smallest; then $G$ has a normal Sylow $p$-subgroup
and a normal Hall $q'$-subgroup.

\begin{cor}
\label{cor:UniqueSgInIntervalOfSS}Let $G$ be a finite supersolvable
group. If $p$ is the smallest prime dividing $[H_{0}:H_{n}]$ for
subgroups $H_{n}\subseteq H_{0}$ of $G$, then there is a unique
subgroup $H_{1}$ with $H_{n}\subseteq H_{1}\subseteq H_{0}$ and
such that $p$ does not divide $\frac{[H_{0}:H_{n}]}{[H_{0}:H_{1}]}$.
\end{cor}
\begin{proof}
Let $\pi=\{ q\,:\, q\leq p,\, q\,|\,\vert H_{0}\vert\}$, and $K$
be a Hall $\pi'$-subgroup of $H_{0}$. Then $K\triangleleft H_{0}$
by the lemma, hence $KH_{n}$ is a subgroup of $H_{0}$ with the desired
properties.
\end{proof}
We are now ready to prove the Main Theorem. The high level idea is
to use the changes in the invariants $I^{p}(H)$ to label cover relations.
Unfortunately, that gives us a lot of identically labelled chains.
So we pick out some distinguished cover relations, and change their
labels to have a unique increasing chain. The details follow.

\begin{thm}
\label{thm:GIscoELshellable}If $G$ is supersolvable with all subgroups
elementary abelian, then $\hat{\mathfrak{C}}(G)$ is co$EL$-shellable,
and so $G$ is coset-shellable. 
\end{thm}
\begin{proof}
We recall by the theorem of P.~Hall restated in Theorem \ref{thm:MainThmRestated}
that $G\subseteq G_{1}\times\dots\times G_{r}$ where each $G_{i}$
is of square free order. If $\hat{\mathfrak{C}}(G_{1}\times\dots\times G_{r})$
is co$EL$-shellable, then it follows immediately from the definition
that the interval $[\emptyset,G]\cong\hat{\mathfrak{C}}(G)$ is as
well. So we can assume without loss of generality that $G=G_{1}\times\dots\times G_{r}$,
the direct product of groups of square free order.

For each $i$, and each $p$ dividing $|G_{i}|$, pick $M_{p,i}^{*}$
to be a maximal subgroup of index $p$ (a Hall $p'$-subgroup) in
$G_{i}$. Such $M_{p,i}^{*}$'s exist because $G_{i}$ is solvable,
and it is a well-known characterization of solvable groups that they
have Hall $p'$-subgroups for each prime $p$ (see \cite[Chapter 9]{Robinson:1996}
for more background). Then set $M_{p,i}$ to be $M_{p,i}^{*}\times\prod_{j\neq i}G_{j}$.
Fix $l(p,j)$ to be an order preserving map into the positive integers
of the lexicographic ordering on the pairs $(p,j)$ for all $p$ dividing
$\vert G\vert$ and $j=1,\dots,r$. We will use $M_{p,j}$ to pick
out the distinguished edges mentioned above. Most edges will be labeled
with $l(p,j)$ for an appropriate $p$ and $j$, while these distinguished
edges will be labeled with the negative of $l(p,j)$.

More precisely, suppose $H_{1}x\subseteq H_{0}x$ is a cover relation
in $\hat{\mathfrak{C}}(G)$. Since $G$ is supersolvable, $[H_{0}:H_{1}]=p$
for some prime $p$, hence the Sylow $p$-subgroups of $H_{1}$ have
dimension (as vector spaces over $\mathbb{Z}_{p}$) one lower than
those of $H_{0}$. It follows that $I^{p}(H_{0})=I^{p}(H_{1})\cup\{ j\}$
for some $j$. Then label the edge $H_{0}x\rightarrow H_{1}x$ as\[
\lambda(H_{0}x\rightarrow H_{1}x)=\left\{ \begin{array}{cl}
-l(p,j) & \textrm{ if }H_{1}x=H_{0}x\cap M_{p,j}\\
\,\,\,\, l(p,j) & \textrm{ otherwise}.\end{array}\right.\]
 Finally, label $\lambda(x\rightarrow\emptyset)=0$. We will show
that $\lambda$ is a co$EL$-labeling.

Intervals in $\hat{\mathfrak{C}}(G)$ all have either the form $[\emptyset,H_{0}x]$,
or $[H_{n}x,H_{0}x]$. We consider these types of intervals separately,
and show there is a unique increasing chain which is lexicographically
first. 

On $[\emptyset,H_{0}x]$, we notice from Proposition \ref{prop:PropertiesOfI}
that every maximal chain from $H_{0}x$ down to $\emptyset$ has 0
on the last edge, and $\pm l(p,j)$ (over all pairs $(p,j)$ such
that $j\in I^{p}(H_{0})$) on the preceding edges. In fact, for each
such pair $(p,j)$, exactly one of $+l(p,j)$ or $-l(p,j)$ occurs
exactly once on any maximal chain. Finally, since $j\in I^{p}(H_{0})$,
there is an element of order $p$ in $H_{0}$ of the form $e_{j}e_{j+1}\cdots e_{r}$,
where each $e_{i}\in G_{i}$ and $e_{j}\neq0$. As $G_{i}\subseteq M_{p,j}$
for $i\neq j$, we see that $H_{0}M_{p,j}=H_{0}M_{p,j}M_{p,j}\subseteq\left(\langle e_{i}\rangle\prod_{i\neq j}G_{i}\right)M_{p,j}=\langle e_{i}\rangle M_{p,j}=G$. 

Then since 0 is the last edge, the only possible increasing chain
is the one with labels $-l(p,j)$ in increasing order. By Lemma \ref{lem:MaxlSubgpsIntersect}
there is such a chain, it is clearly unique and lexicographically
first.

For $[H_{n}x,H_{0}x]$, the situation is only slightly more complicated.
Let a pair $(p,j)$ be called \emph{admissible} for the given interval
if $p$ divides $[H_{0}:H_{n}]$ and $j\in I^{p}(H_{0})\setminus I^{p}(H_{n})$.
If $l(p,j)$ is minimal among admissible $(p,j)$, then there is a
unique $H_{1}x$ of index $p$ in $H_{0}x$ with $H_{0}x\rightarrow H_{1}x$
labeled $\pm l(p,j)$ by Corollary \ref{cor:UniqueSgInIntervalOfSS}
and Lemma \ref{lem:TechnicalLemma}. Moreover, any chain on $[H_{n}x,H_{0}x]$
has exactly one edge with label $\pm l(p,j)$ for each admissible
$(p,j)$.

Suppose $C$ is an increasing chain on $[H_{n}x,H_{0}x]$. Suppose
$H_{i}x\rightarrow H_{i+1}x$ in $C$ is labelled $+l(p,j)$. Then
$l(p,j)$ is minimal among $(p,j)$ admissible for $[H_{n}x,H_{i}x]$
since the chain is increasing. Thus $H_{i}x\rightarrow H_{i+1}x$
is the unique edge down from $H_{i}x$ labeled with $\pm l(p,j)$,
and since the label was positive we see that $H_{n}x\not\subseteq M_{p,j}$.
It follows that the unique increasing chain on $[H_{n}x,H_{0}x]$
is the lexicographically first one labeled with $-l(p,j)$ in increasing
order for $(p,j)$ such that $H_{n}x\subseteq M_{p,j}$, followed
by $+l(p,j)$ for all other admissible $(p,j)$.
\end{proof}

\section{Examples}

At first glance, the labeling constructed in Theorem \ref{thm:GIscoELshellable}
might seem to come {}``from left field.'' It is helpful to work
out what happens for the case where $G$ is a group of square free
order. In this case, many of the complications we faced in the proof
disappear. For example, we don't have to worry about $I^{p}$, since
if $[H_{0}:H_{1}]=p$, then $I^{p}(H_{1})=\emptyset$ and $I^{p}(H_{0})=\{1\}$.
Similarly, we can just take $l(p,j)=p$, since the only possible value
of $j$ is 1. The only $M_{p,j}$'s we have are $M_{p,1}$, which
we can denote as $M_{p}$.

\begin{figure}[htbp]
\includegraphics{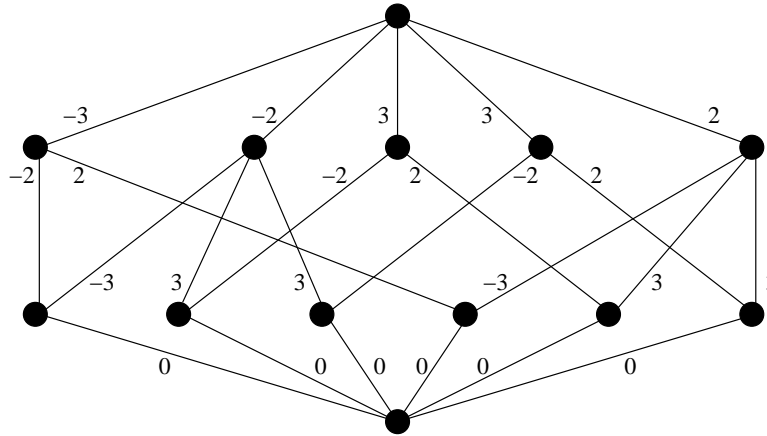}

\caption{The co$EL$-labeling of $\mathfrak{C}(\mathbb{Z}_{6})$. The leftmost
two maximal cosets are $M_{3}=M_{3,1}$ and $M_{2}=M_{2,1}$, respectively.\label{fig:coel_for_z6}}
\end{figure}

Thus we see that for any $H_{0},H_{1}$ with $[H_{0}:H_{1}]=p$ we
get\[
\lambda(H_{0}x\rightarrow H_{1}x)=\left\{ \begin{array}{cl}
-p & \textrm{ if }H_{1}x=H_{0}x\cap M_{p}\\
p & \textrm{ otherwise}\end{array}\right.\]
and $\lambda(x\rightarrow\emptyset)=0$. An example for $\mathbb{Z}_{6}$
is worked out in Figure \ref{fig:coel_for_z6}. An exercise for the
reader might be to work out the labeling for $S_{3}$.

\begin{figure}[htbp]
\includegraphics{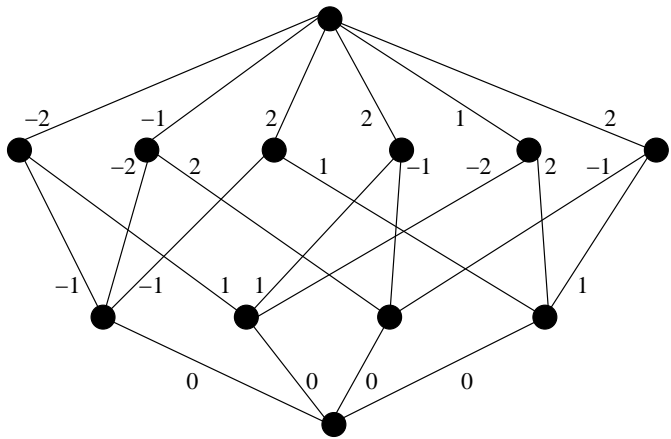}

\caption{The co$EL$-labeling of $\mathfrak{C}(\mathbb{Z}_{2}^{2})$. The
leftmost two maximal cosets are $M_{2,2}$ and $M_{2,1}$, respectively.\label{fig:coel_for_z22}}
\end{figure}

On the opposite extreme, it is not so hard to understand the co$EL$-shelling
on $\mathbb{Z}_{p}^{n}$ -- it is just the change in $I^{p}$, with
$l(p,j)$ becoming $j$. We will not say anything more about this,
but an example for $\mathbb{Z}_{2}^{2}$ is worked out in Figure \ref{fig:coel_for_z22}.%

\section{Consequences and Conclusion}

A (co-)$EL$-labeling of a lattice $L$ tells us a lot about the homotopy
type of $L\setminus\{0,1\}$. In particular, the \emph{falling chains}
(for our purposes, weakly decreasing maximal chains) in an $EL$-labeling
give a basis for the nontrivial homology/cohomology group. See \cite[Section 5]{Bjorner/Wachs:1996}
for a discussion of this in a more general setting. Our co$EL$-labeling
for $\mathfrak{C}(G)$ (where $G$ is a complemented group) thus helps
us understand the cohomology of the order complex in a very concrete
way.

In showing the shellability of a solvable group's subgroup lattice,
Shareshian \cite{Shareshian:2001} produces a so-called {}``coatom
ordering.'' Unfortunately, while the existence of a coatom ordering
implies the existence of something with similar properties to a co$EL$-labeling
(a {}``co$CL$-labeling''), Shareshian is not able to exhibit such
a labeling. Such a labeling would be interesting, as it could presumably
be used to give an alternative proof and/or expand upon a result of
Th{\'e}venaz \cite[Theorem 1.4]{Thevenaz:1985}. Perhaps techniques
like we use here could be used on the chief series for a solvable
group (where every factor is an elementary abelian $p$-group) to
produce a (co-)$EL$-labeling in the subgroup lattice.

\bibliographystyle{/usr/local/teTeX/share/texmf.local/bibtex/bst/hams/hamsplain}
\bibliography{/Users/paranoia/Documents/Research/Master}

\end{document}